\documentclass[reqno,12pt]{amsart}
\usepackage{graphicx}
	\usepackage{color}

\newtheorem{theorem}{Theorem}[section]
\newtheorem{corollary}[theorem]{Corollary}
\newtheorem{lemma}[theorem]{Lemma}
\newtheorem{proposition}[theorem]{Proposition}

\theoremstyle{definition}

\theoremstyle{definition}

\theoremstyle{definition}
\newtheorem{example}[theorem]{Example}

\title{Multiplicity results for bound state solutions of a semilinear equation }
\author{Pilar Herreros}
\address{Departamento de Matem\'atica, Pontificia
        Universidad Cat\'olica de Chile,
       Casilla 306, Correo 22,
        Santiago, Chile.}
\email{\tt pherrero@mat.uc.cl}

\newcommand\R{\mathcal{R}}

\newcommand{\be}{\begin{eqnarray}}
	\newcommand{\ee}{\end{eqnarray}}
\newcommand{\beq}{\begin{equation}}
	\newcommand{\eeq}{\end{equation}}
\newcommand{\ben}{\begin{eqnarray*}}
	\newcommand{\een}{\end{eqnarray*}}

\begin{document}

\begin{abstract}
We consider the problem of multiplicity and uniqueness of radial solutions of a nonlinear elliptic equation of the form
\begin{eqnarray*}
\begin{gathered}
 \Delta u +f(u)=0,\quad x\in \mathbb{R}^N, N\geq 2, \\
 \lim\limits_{|x|\to\infty}u(x)=0.
\end{gathered}
\end{eqnarray*}
where $f$ is a prescribed function, satisfying appropriate conditions.

 This paper is a review of recent developments on the multiplicity of solutions, with a special emphasis on the effect of the mid-range behavior of the nonlinearity $f$ over the existence and multiplicity of solutions.
\end{abstract}

\maketitle

\section{Introduction}

We consider the problem of multiplicity and uniqueness of radial solutions of a nonlinear elliptic equation of the form
\begin{eqnarray}\label{pde}
\begin{gathered}
 \Delta u +f(u)=0,\quad x\in \mathbb{R}^N, N\geq 2, \\
 \lim\limits_{|x|\to\infty}u(x)=0.
\end{gathered}
\end{eqnarray}
where $f$ is a continuous function satisfying appropriate conditions.

Any nonconstant solution to problem \eqref{pde} is called a {\it bound state solution}; nonnegative nontrivial bound states are usually called ground state solutions. Sometimes we use the term {\it higher bound state} to indicate a solution which is not a ground state, or the term $k$-th bound state to refer to a bound state having exactly $k-1$ nodes.

The general problem is to find conditions on this nonlinearity $f$ so that the problem above has either a unique or multiple radial bound state solutions $u: \mathbb{R}^N\to \mathbb{R}$ with a prescribed number of nodal regions. The knowledge of the bound state solutions is very relevant and sometimes crucial in several problems in physics. For example,
it appears in astrophysics as the nonlinear scalar field equation; in fluid mechanics describing
the blowup set of some porous-medium equations, and also in some models in
plasma physics. In addition, if one can determine uniqueness or exact multiplicity  of  $k$-th bound state solutions, it means one has reached full knowledge of the solvability of the problem.

A main step in the study of uniqueness of ground states was given in
the classical work of Gidas, Ni and Nirenberg'79 \cite{gnn1, gnn2}. They show
that all ground states of this problem are radially symmetric for instance when $f\in C^{1+\mu}$, $\mu>0$, $f(0)=0$ and $f'(0)<0$. Since then, several works in this direction have been published. See, for instance, the paper of Li and Ni'93 \cite{ln} and the references therein.

In view of these results it is natural to focus our attention on radially symmetric solutions,
and the radial version of \eqref{pde}, that is
\begin{eqnarray}\label{eq2}
\begin{gathered}
u''+\frac{N-1}{r}u'+f(u)=0,\quad r>0,\quad N\geq 2,\\
u'(0)=0,\quad \lim\limits_{r\to\infty}u(r)=0,
\end{gathered}
\end{eqnarray}
where  $'$ denotes differentiation with respect to $r$.

The goal of this article is to show how the different properties of the nonlinearity $f$ affect the existence and uniqueness of bound state solutions. We will begin discussing some of the classical results and tools that have been used to study this problem, where the conditions imposed on $f$ relate mainly to the behavior of $f$ near $0$, and different growth assumptions. Our main objective is to review two recent articles where the type of conditions imposed on $f$ are different, relating to the behavior of $f$ in its middle range.  To this end we will consider functions $f$ that either change sign several times or increase fast in the middle.

The article is organized as follows. In Section $2$ we discuss some of the classical approaches to this problem, trying to give an idea of the importance of the conditions used and a description of some tools used in these approaches. This does not intend to be a comprehensive review of all results, not even of the most important ones. The choice of which results and techniques are included is based on the goal of understanding the last sections.
This Section is divided into three subsections, the first one focusing on existence of ground state solutions, where the conditions on $f:[0,\gamma)\to \R$ near $0$ and near $\gamma$ are given, where $0<\gamma\leq \infty$. The second subsection is about the uniqueness of solutions, where growth assumptions are needed, mentioning one approach for its proof as well as some cases where there is multiplicity of solutions. The third subsection gives an idea of how to extend some of the previous results to $k$th-bound state solutions.

Section $3$ is devoted to the case where $f$ changes sign several times, in such a way that its primitive $F$ has at least one positive local maximum in the interior. By studying the solutions with initial condition $u(0)$ bigger than this maximum, we can see that there must be at least $2$ such $k$th-bound state solutions for each $k>k_0$, for some $k_0\geq 1$. On the other hand, ground state solutions $(k=1)$ might not exist, as a lower bound on $k_0$ can be given in terms of the local maxima and minima of $F$.

Section $4$ deals with the case where $f$ changes sign only once, but has an abrupt magnitude change. We consider $f$ defined by parts, with nice regularity and growth assumptions for small values of $u$, up to a height $\alpha_*$ that corresponds to a ground state solution. Above this point, we make an arbitrarily high jump and continue with any (large) function, glued by a small linear interpolation for continuity.
We can see that this jump generates a second ground state solution. Moreover, if we follow this jump with another one, this time making $f$ small, we can find a third ground state solution (and thus as many as we want).

\section{Some Background and past approaches}

\subsection{Existence of ground state solutions}

Consider the initial value problem
\begin{eqnarray}\label{ivp}
\begin{gathered}
u''+\frac{N-1}{r}u'+f(u)=0,\quad r>0,\quad N\ge 2,\\
u(0)=\alpha ,\quad u'(0)=0
\end{gathered}
\end{eqnarray}
We will assume $f$ is an odd function, defining $f$ only for $s\geq0$ and extending it by $f(-s)= - f(s)$. Thus when we say, for example, $f(s)=s^p$ we mean $f(s)=|s|^{p-1}s$.

The first property we can see of these solutions is that for critical points $u''=-f(u)$, thus it will have local maxima when $f>0$ and minima when $f<0$. Also, if when $r\to \infty$  we have
$$ \lim_{r\to \infty} u(r)= L  , \quad \lim_{r\to \infty} u'(r)= 0 \quad \mbox{ and} \quad  \quad \lim_{r\to \infty} u''(r)= 0 $$
then  $f(L)=0$. Thus it only makes sense to look for ground states if $f(0)=0$.

When $f(s)>0$ for $s>0$, solutions to the initial value problem will begin at a maximum in $u(0)=\alpha$ and either stay positive and tend to $0$ at infinity, thus being a ground state, or will change sign at some $r_0>0$ where $u(r_0)=0$. Which of these behaviors occur will depend on the growth of $f$. It is well known, for instance, that for $f(s)=s^p$ there is a critical exponent $p^*=\frac{N+2}{N-2}$ that separate these cases: If {$p<p^*$} there are no ground state solutions while if $p>p^*$ all initial conditions give a ground state solution.

One of the classical ways of proving existence of these solutions is the use of Pohozaev's Identity \cite{poh}
$$r^Nu'^2 + (N-2)r^{N-1}u'u +2r^N F(u)= \int_0^r r^{N-1}\left( 2NF(u)-(N-2)f(u)u \right) \ dt $$
where  $\displaystyle F(u)=\int_0^uf(t)\ dt$. If there is a point where $u(r)=0$ it gives
$$0 \leq  r^Nu'^2 = \int_0^r t^{N-1}\left( 2NF(u)-(N-2)uf(u) \right) \ dt  = \int_0^r t^{N-1}Q(u(t)) \ dt $$
Thus if {$Q(s)= 2NF(s)-(N-2)sf(s) <0$} solutions are {positive}. For $f(s)=s^p$, $Q(s)= \left(\frac{2N}{p+1}-(N-2)\right)s^{p+1}$ will be negative when {$p>p^*=\frac{N+2}{N-2}$}, and solutions will be ground states.

A more interesting behavior can be seen when $f$ changes sign. If we consider functions $f$ that are negative for $0<s<b$ and positive for $s>b$, a solution to the initial value problem with $u(0)=\alpha>b$ will start at a maximum, but it can have a local minimum when $u<b$. The best studied function with this shape is $f(s)=s^p-s$ with $p>1$. If we want to use Pohozaev's identity we have $Q(s)= \left(\frac{2N}{p+1}-(N-2)\right)s^{p+1}-2s^2$ that will be negative when {$p>p^*=\frac{N+2}{N-2}$} or when {$s$ is small}. Unlike the previous case, the fact that $u(r)>0$ does not imply that they are ground states, on the contrary they will have a local minimum and oscillate around $1$, which is the other point with $f(1)=0$.

\begin{figure}[h]
\begin{center}
 \includegraphics[keepaspectratio, width=6cm]{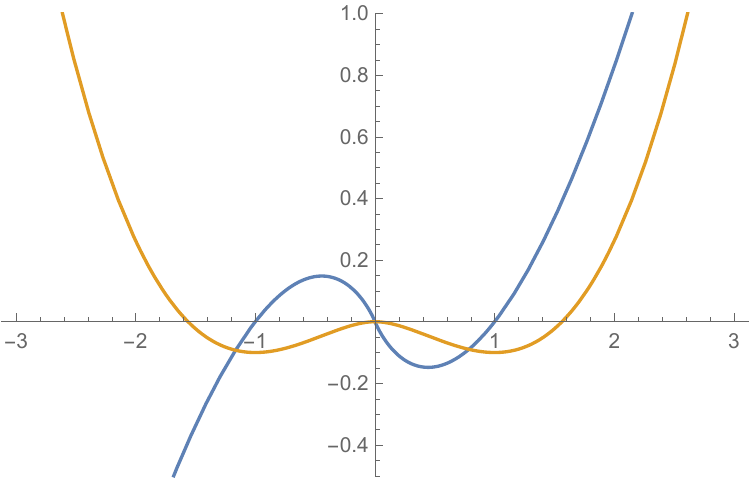} \includegraphics[keepaspectratio, width=6cm]{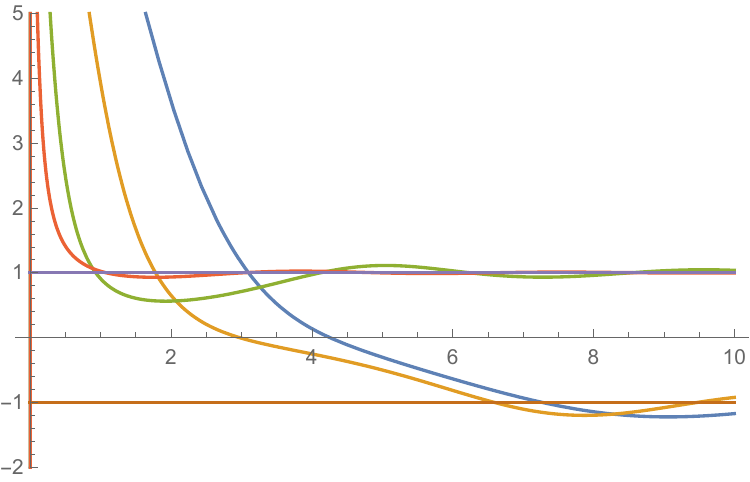}
 \caption{${f(u)}=u^p-u$ and {$F(u)$} \hspace{0.5cm} Solutions for $N=4$, $\alpha=10$, different values of $p$}\label{Fig-up-u}
 \end{center}
 \end{figure}

\begin{figure}[h]
\begin{center}
 \includegraphics[keepaspectratio, width=10cm]{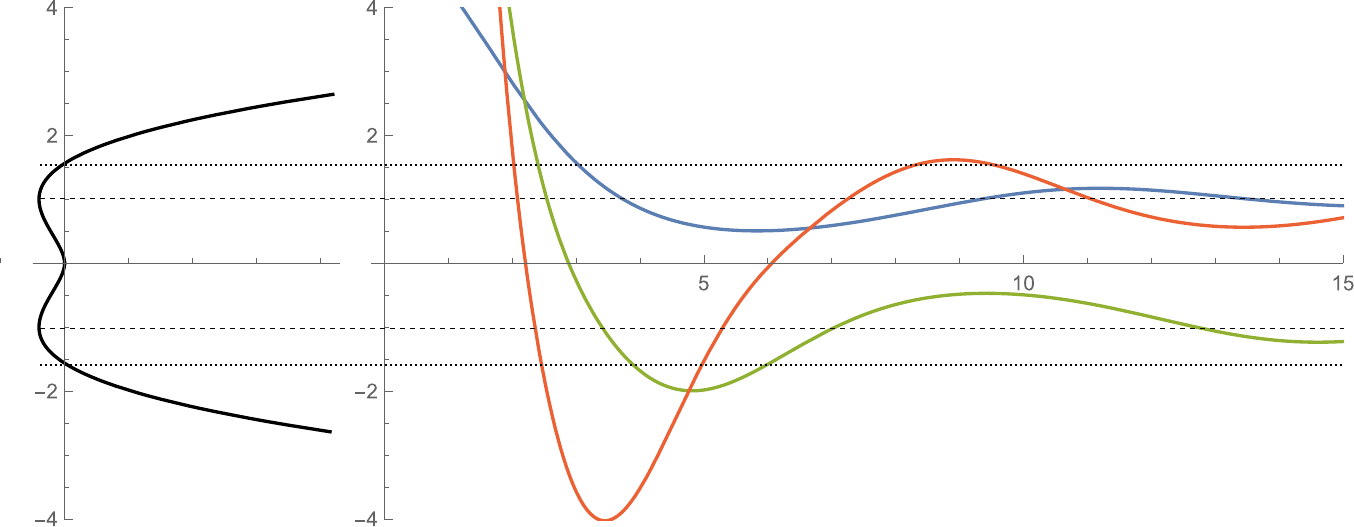}
 \caption{ $F(u)$ and solutions for $N=4$, $f(s)=s^{\frac{3}{2}}-s$, different values of $u(0)=\alpha$.}\label{Fig-FySolup-u}
  \end{center}
 \end{figure}

 The biggest difference occurs in the subcritical case, when $p<p^*$, that the behavior of solutions depends on the initial condition $u(0)=\alpha$. When $\alpha$ is small solutions will be positive, while when $\alpha$ is big they will change sign (See Figure \ref{Fig-up-u}). A useful tool for this case is the energy functional
\begin{eqnarray*}
{I(r)=\frac{|u'(r)|^2}{2}+F(u(r))} \quad \mbox{with} \quad {I'(r)=-(N-1)\frac{|u'(r)|^2}{r}\leq 0}.
\end{eqnarray*}
If we think of $r$ as time, $I$ can be thought of as having a potential part given by the ``height'' of $F(u)$ and a kinetic part depending on speed. The total energy is decreasing, and decreases more slowly for larger $r$. If a solution reaches $0$ at an $r_0>0$ we have
$$F(u(0))=I(0)>I(r_0)=\frac{|u'(r_0)|^2}{2}>0,$$ so solutions with initial condition $u(0)=\alpha \leq\beta$ will be positive, where $\beta=\left(\frac{p+1}{2}\right)^{\frac{1}{p-1}}$ is where $F(\beta)=0$. Note also that if a solution has a positive local minimum at $r_0$, then $F(u(r_0))<0$ and the solution can not reach $0$ after that. (See Figure \ref{Fig-FySolup-u})

This argument depends only on the general shape of $f$ near $0$, we will consider from now on functions $f$ with this general shape by imposing that it satisfies the following condition.
\begin{itemize}
\item[$(C_1)$] $i)$ $f\in C([0,\gamma])\cap C^1((0,\gamma])$ with  $f(0)=0$, $f(s)<0$ for $0<s<\epsilon$,  \\
$ii)$ $\exists b>0$ such that $f(s)\leq 0$ in $[0,b]$ and $f(s)> 0$ for $s>b$\\
$iii)$ $\exists \beta>0$ such that $F(s)\leq 0$ in $[0,\beta]$ and $F(s)> 0$ for $s>\beta$.
\end{itemize}
where $\gamma\in (0,\infty]$  considering $[-\infty,\infty]$ as $\R$.

Let us denote by $u(r,\alpha)$ the solution of \eqref{ivp} with $u(0,\alpha)=\alpha$. By the above argument solutions that have a positive local minimum will stay positive after that, while other solutions may exist that change sign. We will classify the initial conditions by the behavior of the corresponding solutions in the following way.
Let us set
$$Z_1(\alpha):=\sup\{r>0\ |\ u(s,\alpha)>0\mbox{ and }u'(s,\alpha)<0\ \mbox{ for all }s\in(0,r)\}$$
and define
\begin{eqnarray*}
{\mathcal P_1}&=&\{\alpha\in[\beta,\gamma)\ :\ u(Z_1(\alpha),\alpha)>0\}\\
{\mathcal N_1}&=&\{\alpha\in[\beta,\gamma)\ :\ u(Z_1(\alpha),\alpha)=0\quad\mbox{and}\quad u'(Z_1(\alpha),\alpha)<0\}\\
{\mathcal G_1}&=&\{\alpha\in[\beta,\gamma)\ :\ u(Z_1(\alpha),\alpha)=0\quad\mbox{and}\quad u'(Z_1(\alpha),\alpha)=0\}\\
\end{eqnarray*}
We note that, if $f$ is differentiable at $0$, or Lipschitz, the only solution with $u(r_0)=u'(r_0)=0$ is the constat $u\equiv 0$. Therefore if $\alpha^*\in {\mathcal G_1}$ we must have $Z_1(\alpha^*)=\infty$ and $u(r,\alpha^*)$ is a ground state. If not, at a point with $u(r_0)=u'(r_0)=0$ the solution will stop being unique, and we can continue the solution as $0$ afterwards. In this case we will also consider the solution as a ground state, sometimes called  ground state with compact support.

Since the solutions to the initial value problem depend continuously on $\alpha$, the sets ${\mathcal N_1}$ and ${\mathcal P_1}$ are open sets.  Also, $\beta\in {\mathcal P_1}$ so  ${\mathcal P_1}$ is nonempty. To prove existence of solutions, it is enough to prove that ${\mathcal N_1}$ in nonempty, so that  ${\mathcal N_1}$ and ${\mathcal P_1}$ cannot cover $[\beta,\gamma]$ and thus there is an $\alpha^*\in{\mathcal G_1}$.

{Existence} of ground state solutions can be obtained under different ``subcritical type conditions''. See for example: Atkinson and Peletier'$86$ \cite{ap}, Berestycki and Lions'$83$ \cite{b-l1}, De Figueiredo and  Ruf'$95$ \cite{dfr}, Ferrero and Gazzola \cite{fg}.
The most general one, up to our knowledge, is the work of F. Gazzola, J. Serrin and M. Tang'$00$ \cite{gst} using a condition that was considered first by Castro and Kurepa'87 \cite{cku}. They ask only for a subcritical condition near $\gamma$, obtaining existence under the following condition
\begin{itemize}
\item[(GST)]    If $\gamma<\infty$ then $f(\gamma)=0$ or \\
 if $\gamma=\infty$, $Q(s)$ is locally bounded below near $s = 0$ and there exists $\bar b>\beta$  and $k\in(0,1)$ such that
 $Q(s)>0$ for all $s>\bar b$,
\ben\limsup_{s\to\infty}Q(s_2)\Bigl(\frac{s}{f(s_1)}\Bigr)^{N/2}=\infty,
\een
for all $s_1, s_2\in [ks,s]$, where $Q(s):=2NF(s)-(N-2)sf(s)$.\\
\end{itemize}

Note that this condition for $\gamma=\infty$ is satisfied for $s^p$ with  $p<p^*$, and any function that behaves like it near infinity.

The main ingredient of this proof is the idea that, since the energy decreases more slowly when $r$ is big, a solution that reaches a certain height too far will have too much energy to stop before reaching $0$. This is formalized in the following lemma that was proved for more general operators, under general shape conditions that are weaker than $(C_1)$ (See \cite[Lemma 3.1]{gst}).

\begin{lemma}[Gazzola-Serrin-Tang]\label{LemmaGST}
Given $s_0>\beta$, let $r_0$ be where $u(r_0)=s_0$. If
$${r_0\geq C(s_0)}=\sqrt{2}(N-1)\frac{s_0}{F(s_0)}\left(F(s_0)-\min_{0<s<\beta}F(s) \right)^{\frac{1}{2}} $$
then the {solution changes sign}.\\
\end{lemma}

They conclude the proof of existence showing that condition $(GST)$, among others for more general equations, is enough to prove that as $s_0\to \gamma$ there will be solutions with $r_0\geq C(s_0)$. Thus ${\mathcal N_1}\neq \emptyset$ and there is a ground state solution with initial condition between ${\mathcal N_1}$ and ${\mathcal P_1}$.

\subsection{Uniqueness and multiplicity of ground states}

{Uniqueness} has been established for functions of general shape as in condition $(C_1)$ with some growth assumption, usually  superlinear or sublinear growth. Key steps were given in the work of Coffman'72 \cite{coff}, Peletier and Serrin'83 \cite{pel-ser1}, Kwong'89 \cite{Kw} and
Franchi, Lanconelli and Serrin'96 \cite{fls}. See also McLeod and Serrin'87 \cite{ms}, Chen and Lin'91 \cite{cl}, Pucci and Serrin'98 \cite{pu-ser}, Cort\'azar Elgeta and Felmer'98 \cite{cfe2}, Serrin and Tang'00 \cite{st}.

Several approaches have been used to compare any two solutions to the initial value problem,  proving either that they do not intersect, or that they intersect at most once, in certain intervals and thus are ``ordered'' in some way.

We will mention one approach in particular, used by  Erbe and Tang in \cite{et} to prove uniqueness of solutions for $f>0$ and superlinear, satisfying  $$\left(\frac{F}{f}\right)'\geq \frac{N-2}{2N} \quad \mbox{ and } \quad f(u)<{u f'(u)}. $$
 Since solutions are monotonously decreasing in $[0,Z_1(\alpha)]$, we can study instead the inverse $r(s)$ for $s\in[u(Z_1(\alpha), \alpha]$ that satisfies the equation
$$r''(s)-\frac{N-1}{r(s)}(r'(s))^2 - f(s)(r'(s))^3=0.$$
This approach facilitates the comparison between solutions since the value of $f(s)$ does not depend on the solution.

 Erbe and Tang prove uniqueness using a Pohozaev type functional
\begin{equation}\label{PdeET}P(s)=-2N\frac{F}{f}(s)\frac{r^{N-1}(s)}{r'(s)}-\frac{r^N(s)}{(r'(s))^2}
-2r^N(s)F(s),\quad \end{equation}
with
\begin{equation*}P'(s)=\frac{\partial P}{\partial s}(s)=\left(N-2-2N\Bigl(\frac{F}{f}\Bigr)'(s)\right)\frac{r^{N-1}(s)}{r'(s)},
\end{equation*}
they compare any two solutions to the initial value problem, proving that they can intersect only once before reaching $0$.

On the opposite side, there are nonliniarities $f$ such that problem \eqref{eq2} does not have a unique solution. The {multiplicity}  problem has been studied mostly for the  non-autonomous case $f(x,u)=g(x,u)-a(x)u$,
for different nonnegative functions $g$ and coefficients $a$.
See for example: Cao and Zhou '96 \cite{cao}, Adachi and Tanaka'00 \cite{at1}, Cerami, Hsu and Lin'10 \cite{hl}, Wei and Yan'10 \cite{wy},  Ao and Wei'14 \cite{aw}, Del Pino, Wei and Yao'15 \cite{dwy}, Cerami and Molle '19 \cite{cm19}, Molle and Passaseo'21 \cite{mp21}.

Some progress has also been made for the autonomous case, for example  D\'avila, del Pino and Guerra'13 \cite{ddg} proved that for {$f(u)=-u+u^p+\lambda u^q$} with $N=3$, $1<q<3$, $p<5$ near $5$, if $\lambda$ is large enough, then there exist {at least three} radial ground state solutions to this problem.  Wei and Wu'22 \cite{wei-wu}, considered the nonlinearity {$f(u)=|u|^{2^*-2}u+\lambda u+\mu |u|^{q-2}u$}, where $2^*=\frac{2N}{N-2}$, $N=3$, $2<q<10/3$, and proved that under some conditions in $\mu>0$ the problem has a {second ground state} for some $\lambda<0$.

\subsection{Existence and uniqueness of $k$th-bound states}

\begin{figure}[h]
\begin{center}
 \includegraphics[keepaspectratio, width=10cm]{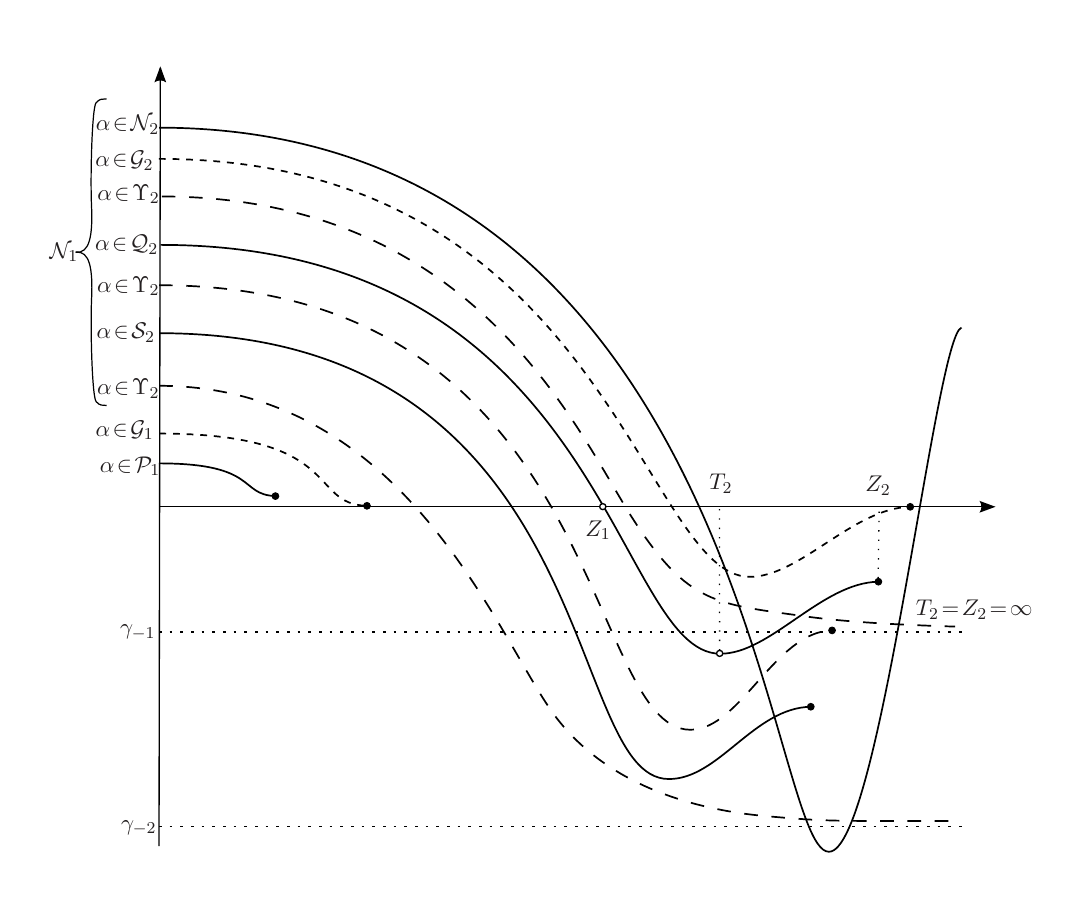}
 \end{center}
 \caption{Solutions of \eqref{ivp} with initial condition in these sets.}\label{Fig-Conj2}
\end{figure}

To extend some of the existence results to $k$th-bound state solutions, we will continue the analysis of a  solution to the initial value problem after it crosses the value $0$. These solutions can either decrease forever, in which case it must tend to the point $b$ where $f(b)=0$, or have a local minimum. We will extend by induction the sets $\mathcal  N_1$, $\mathcal P_1$ and $\mathcal G_1$  in the following way. (See Figure \ref{Fig-Conj2})

If ${\mathcal N_{k-1}}\not=\emptyset$, we set
$${\mathcal F}_k=\{\alpha\in\mathcal N_{k-1}\ :\ (-1)^ku'(r,\alpha)\le 0\quad\mbox{for all }r> Z_{k-1}(\alpha)\}.$$
For $\alpha\in \mathcal N_{k-1}\setminus{\mathcal F}_k$, we set
\ben
T_{k-1}(\alpha):&=&\sup\{r\in(Z_{k-1}(\alpha),D_\alpha)\ :\ (-1)^ku'(r,\alpha)\le 0\},
\een
and for $\alpha\in{\mathcal F}_k$, we set $T_{k-1}(\alpha)=\infty$.

Next, for $\alpha\in \mathcal N_{k-1}\setminus {\mathcal F}_k$, we define the extended real number
\begin{eqnarray*}
Z_k(\alpha):=\sup\{r>T_{k-1}(\alpha)\ |\ (-1)^ku(s,\alpha)<0\mbox{ and }(-1)^ku'(s,\alpha)>0\ \\
\mbox{ for all }s\in(T_{k-1}(\alpha),r)\},
\end{eqnarray*}
and again if $\alpha\in{\mathcal F}_k$, we set $Z_k(\alpha)=\infty$.

We now define
\begin{eqnarray*}
{\mathcal N_k}&=&\{\alpha\in\mathcal N_{k-1}\setminus{\mathcal F}_k\ :\ u(Z_k(\alpha),\alpha)=0\quad\mbox{and}\quad (-1)^ku'(Z_k(\alpha),\alpha)>0\},\\
{\mathcal G_k}&=&\{\alpha\in\mathcal N_{k-1}\setminus {\mathcal F}_k\ :\ u(Z_k(\alpha),\alpha)=0\quad\mbox{and}\quad u'(Z_k(\alpha),\alpha)=0\},\\
{\mathcal P_k}&=&\{\alpha\in\mathcal N_{k-1}\ :\ (-1)^ku(Z_k(\alpha),\alpha)<0\}.
\end{eqnarray*}

Using similar arguments as in \cite{gst}, C. Cort\'azar, M. Garc\'ia-Huidobro and C. Yarur'13  \cite{cghy4} proved that for each $k$ the solutions with large initial condition are in $\mathcal N_k$, therefore we can define an $\alpha_k$ as the minimum $\alpha$ such that $(\alpha,\gamma_*)\in\mathcal  N_k$. To conclude that this solution is a bound state in $\mathcal G_k$ we also need to prove that $\mathcal P_k$ is non-empty and $\alpha_k\neq \alpha_{k-1}$. This can be achieved by induction, showing that near an $\alpha_k\in \mathcal G_k$ solutions won't have enough energy to cross $0$ another time, so there will be a neighborhood $(\alpha_k-\epsilon,\alpha_k+\epsilon)$ of elements not in $\mathcal N_{k+1}$. Since $(\alpha_k,\alpha_k+\epsilon)\subset(\alpha,\gamma_*)\in\mathcal N_k\setminus \mathcal N_{k+1}$, there must be elements in $\mathcal P_{k+1}$ and $\alpha_{k+1}>\alpha_k$.

Some of the uniqueness results of ground states can also be extended to bound state solutions.  C. Cort\'azar, M. Garc\'ia-Huidobro and  C. Yarur \cite{cghy,cghy2} proved uniqueness of  $2$th-bound  state solutions when $f$ satisfies $ f(s)\leq f'(s)(s-b) $ and
$$ \frac{s f'(s)}{f(s)} \quad \mbox{ nonincreasing for}\quad  s>b \quad \mbox{and } \frac{\beta f'(\beta)}{f(\beta)}\leq\frac{N}{N-2},$$
and uniqueness of  $k$th-bound  state solutions under stronger assumptions.

\section{Multiplicity via sign changes of $f$}

In this section we want to see how some zeros of $f$, or more precisely some local maxima of $F$, can generate multiplicity of solutions. This was studied by the author in collaboration with  Carmen Cort\'azar and Marta Garc\'ia Huidobro in \cite{cghh15}, and we refer to this article for more details.

Consider functions $f$ that satisfy
\begin{enumerate}

\item[$(A_1)$] $f$ is a continuous odd function defined in $[-\gamma_*,\gamma_*]$,  $f(0)=0$, and $f$ is locally Lipschitz in $[-\gamma_*,\gamma_*]\setminus\{0\}$.
\item[$(A_2)$]
There exists $\delta>0$ such that if we set $F(s)=\int_0^sf(t)dt$, it holds that $F(s)<0$ for all $0<|s|<\delta$, and $F(s)<F(\gamma_*)$ for all $s\in(-\gamma_*,\gamma_*)$.
\item[$(A_3)$] $F$ has a local maximum at some $\gamma\in(\delta,\gamma_*)$ with $F(\gamma)>0$.
\item[$(A_4)$] $f$ has a finite number of zeros in $(\delta,\gamma_*)$ and $f$ changes sign at these  points.
\item[$(A_5)$] $\gamma_*<\infty$ with $f(\gamma_*)=0$ or $\gamma_*=\infty$ with
$$\lim_{s\to\infty}\Bigl(\inf_{s_1,s_2\in[\theta s,s]}Q(s_2)\Bigl(\frac{s}{f(s_1)}\Bigr)^{N/2}\Bigr)=\infty,$$ where $Q(s):=2NF(s)-(N-2)sf(s)$.

\end{enumerate}

\begin{figure}[h]
\begin{center}
 \includegraphics[keepaspectratio, width=8cm]{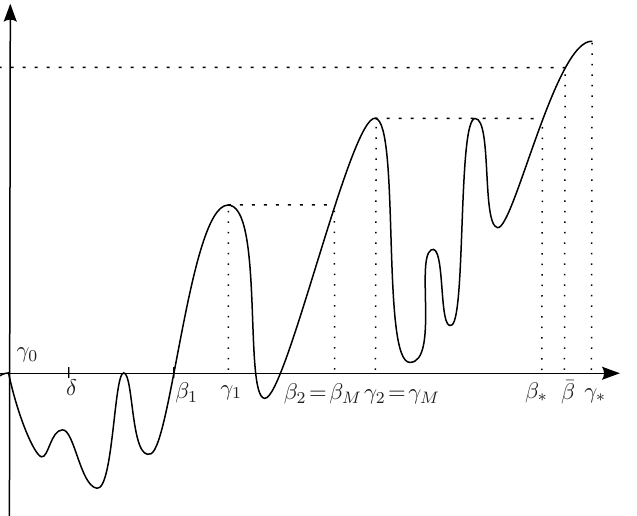}
 \end{center}
  \caption{$F(u)$ and the constants $\gamma_i$ and $\beta_i$.}\label{Fig-defjb}
\end{figure}
We note that $(A5)$ is a "subcritical at $\gamma_*$" condition equivalent to condition $(GST)$, and for $\gamma_*=\infty$ it is satisfied by $s^p$ with $p<p^*$.

Since $f$ satisfies $(A_4)$, its primitive $F$ has a finite number of critical points and they are all local maxima or minima. We will enumerate some of the local maxima in the following way. Denote by $\gamma_1$ the first local maximum of $F$ with $F(\gamma_1)>0$, and by $\gamma_i$ first local maximum of $F$ with $F(\gamma_i)>F(\gamma_{i-1})$ for $i=1,2,\dots, M$ (See Figure \ref{Fig-defjb}). Thus $F(\gamma_M)$ is the highest local maximum of $F$ in $(0,\gamma_*)$, lets denote by $\beta_*$ the largest point in $(\gamma_M,\gamma_*)$ where $F(\gamma_M)=F(\beta_*)$. We will study the behavior of solutions with initial value in $(\beta_*,\gamma_*)$.

 The most important difference between this case and the previous one is that solutions with initial condition close to $\beta_*$ will never cross the height $s=\gamma_M$, staying in this range forever. Moreover, solutions with enough initial energy (big enough $\alpha$) to cross the height $s=\gamma_M$ might have too much energy to stop before reaching $0$, so in many cases there will be no ground state solutions, or $k$-bound state solutions for small $k$.

The main result is the existence of two $k$th-bound state solutions with initial value in $(\beta_*,\gamma_*)$ for $k$ large enough.

\begin{theorem}\label{main0}
Let $f$ satisfy conditions $(A_1)-(A_5)$ above, then
there exists $k_0\in\mathbb N\cup\{0\}$ such that for any $k\ge k_0$, there exist at least two solutions $u$ of \eqref{eq2}, with initial value in $(\beta_*,\gamma_*)$, having exactly $k$ sign changes in $(0,\infty)$.

\end{theorem}

Note that for $i>1$ we have $f(\gamma_i)=0$ and  $f$ restricted to $[-\gamma_i, \gamma_i] $ satisfies all the hypothesis of Theorem \ref{main0}. Therefore we will have  two solutions $u$ of \eqref{eq2}, with initial value in $(\beta_i,\gamma_i)$, having exactly $k$ sign changes for $k>k_i$. Furthermore, by the classical results seen in Section $2$, there will be one solution with initial value in $(\beta_1,\gamma_1)$ for all $k$.  This result allows us to give conditions on $f$ so that problem \eqref{pde} has at least any given number of solutions having a exactly $k\ge k_0=\max k_i$ nodes.\\

\begin{corollary}
Assume that $f$ satisfies assumptions $(A1)-(A5)$. Then there exists $k_0\in\mathbb N\cup\{0\}$ such that for any $k\ge k_0$, there exist at least $2M+1$ solutions of \eqref{eq2} with a positive initial value having exactly $k$ sign changes in $(0,\infty)$.
\end{corollary}

We also prove that $k$th-bound states do not exist in general for small $k$, giving a lower bound on $k_0$.

\begin{theorem}
 If $f$ satisfies assumptions $(A1)-(A5)$ and
\beq\label{t13}-\min_{s\in[0,\beta_{*}]}F(s)<\frac{(\beta_*-\gamma_1)}{2(N-1)(k+1)}\frac{F(\gamma_1)}{2\gamma_1}-F(\gamma_*),
\eeq
then there are no solutions $u$ of \eqref{eq2}, with initial value in $(\beta_*,\gamma_*)$, having exactly $j$ sign changes in $(0,\infty)$ for any $j=0,1,\ldots, k$.

When $\gamma_*=\infty$, $F(\gamma_*)$ should be replaced by $\sup_{s\in[0,\alpha_k]}F(s)$ where $\alpha_k$ is such that $[\alpha_k,\gamma_*)\in \mathcal N_k$.
\end{theorem}

On the other hand, we give a sufficient condition on $f$ so that $k_0=1$ in Theorem \ref{main0} above, showing the existence of ground states. Unfortunately, these are complicated conditions as seen below.

\begin{theorem} If
\beq\label{caso1}
(\bar C+A)\bar I
<\frac{2^{1/2}(N-1)}{(\bar I+\bar F)^{1/2}}\int_0^{\beta_1} |F(s)|ds,
\eeq
then for {all $k\geq 0$} there exist {at least two solutions} $u$ of \eqref{eq2}, with initial value in $(\beta_*,\gamma_*)$, having exactly $k$ sign changes in $(0,\infty)$, where

$$\bar F:=-\min_{s\in[0,\beta_1]}F(s)>0.$$
{If} $\gamma_*<\infty,$ $$\quad\mbox{}\quad   A=\frac{\beta_*-\beta_1}{((F(\bar\beta)-F(\gamma_M)))^{1/2}}+ \Bigl(\frac{2N(\bar\beta-\beta_*)}{\min_{t\in[\beta_*,\bar\beta]}f(t)}\Bigr)^{1/2}\quad\mbox{and}\quad \bar I=F(\gamma_*), $$
{if} $\gamma_*=\infty,$ $$\quad\mbox{}\quad   A=\max\{1, \frac{\beta_*-\beta_1}{((F(\bar\beta)-F(\gamma_M)))^{1/2}}+ \Bigl(\frac{2N(\bar\beta-\beta_*)}{\min_{t\in[\beta_*,\bar\beta]}f(t)}\Bigr)^{1/2}\}\quad\mbox{and}\quad $$
\begin{eqnarray*}
\bar I:=\Bigl(\frac{\bar C+1}{\bar C}\Bigr)^N\Bigl(2F(\bar\beta)+(\bar\beta-\beta_1)^2+\frac{1}{N}\Bigl(\sup_{s\in[\beta_1,\bar\beta]}Q(s)
-\min_{s\in[s_0,\bar\beta]}Q(s)\Bigr)\Bigr)\\
+\frac{(N-2)^2\bar\beta^2}{2\bar C^2}\end{eqnarray*}
 $$\quad\mbox{ with}\quad  \bar C:=2(N-1)\frac{\bar\beta-\beta_1}{F(\bar\beta)-F(\gamma_M)}(2(F(\bar\beta)-\min_{s\in[\beta_1,\beta_{*}]}F(s)))^{1/2}.$$

\end{theorem}

The proof of Theorem \ref{main0} follows some ideas mentioned in the previous section, but we have to deal with the presence of local maxima of $F$. This local maxima will act as obstacles, being determinant in the behavior of solutions to the initial value problem and the multiplicity of bound states.

\begin{figure}[h]
\begin{center}
 \includegraphics[keepaspectratio, width=12cm]{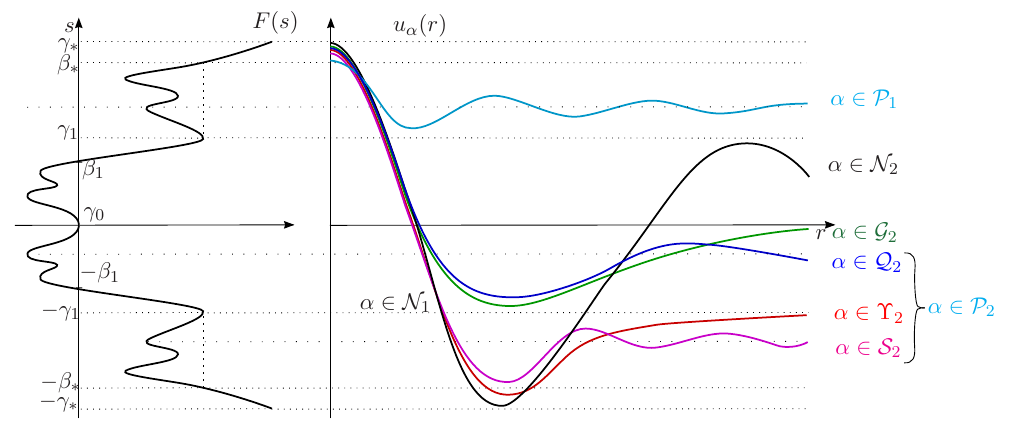}
 \end{center}
 \caption{$F(u)$ and solutions for different values of $u(0)$}\label{Fig-FySol2}
\end{figure}

To study these solutions in \cite{cghh15} we separate positive solutions depending on where they end, we decompose each set $\mathcal P_k$ into (See Figure \ref{Fig-FySol2})
\begin{eqnarray*}
\mathcal Q_k&=&\{\alpha\in\mathcal P_k\ :\  0<|u(Z_k(\alpha),\alpha)|<\gamma_1\}\\
\mathcal S_k&=&\bigcup_{i=1}^{M}\{\alpha\in\mathcal P_k\ :\ \gamma_i<|u(Z_k(\alpha),\alpha)|<\gamma_{i+1}\}\\
\Upsilon_k&=&\bigcup_{i=1}^{M}\{\alpha\in\mathcal P_k\ :\ |u(Z_k(\alpha),\alpha)|=\gamma_i\}
\end{eqnarray*}
where the constants $\gamma_i$ are the increasing sequence of local maxima with $F(\gamma_1)>0$ and $F(\gamma_i)<F(\gamma_{i+1})$.

As the minima (maxima) of $u$ occur at values where $f(u)\le 0$ ($f(u)\ge 0$), it follows that if $\alpha\in \mathcal S_j\cup\mathcal Q_j$, with $\gamma_i<|u(Z_j(\alpha),\alpha)|<\gamma_{i+1}$, then $F(u(Z_j(\alpha),\alpha))<F(\gamma_i)<F(\gamma_{i+1})$ and hence $\gamma_i<|u(r,\alpha)|<\gamma_{i+1}$ for all $r>Z_j(\alpha)$.
It should be noticed that, by the unique solvability of \eqref{ivp} up to a double zero,  if $\alpha\in\Upsilon_k$ then necessarily $Z_k(\alpha)=\infty$. Using a modification of Lemma \ref{LemmaGST} we can prove that near a solution with $\alpha\in \Upsilon_k$, solutions cross $0$ as many times as we want, thus $(\alpha-\epsilon_j,\alpha+\epsilon_j)\subset \mathcal N_j\cup \Upsilon_k$ for all $j$. We can also define as before $\alpha_k$ as the infimum such that $(\alpha_k,\gamma*)\in  \mathcal N_k$.
If $\mathcal G_k\neq\emptyset$ we can define our candidates as
$${\underline{\alpha}}_{k}:=\inf(\mathcal G_{k}\cup\mathcal Q_{k})\quad \mbox{and}\quad\bar\alpha_{k}:=\sup(\mathcal G_{k}\cup\mathcal Q_{k}).$$
By the above argument, and since $\mathcal S_k, \mathcal Q_k$ and $\mathcal N_k$ are open, these elements  must be in $\bigcup_{i=1}^k\mathcal G_i$.

Using similar arguments as in \cite{cghy4} we can see that near $\alpha_k\in\mathcal  G_k$ solutions will not have enough energy to cross $0$ another time nor $\gamma_1$, so there will be a neighborhood such that $(\alpha_k-\epsilon,\alpha_k+\epsilon)\subset \left( \mathcal Q_{k}\cup  \mathcal G_{k}\cup  \mathcal Q_{k+1}\right)$. Thus, just below $\underline{\alpha}_{k}$ and just above $\bar{\alpha}_{k}$ there must be elements of $\mathcal Q_{k+1}$, ensuring that $\underline{\alpha}_{k+1}<\underline{\alpha}_{k}$ and
$\bar{\alpha}_{k}<\bar{\alpha}_{k+1}$.  Moreover, $\underline{\alpha}_{k+1}, \bar{\alpha}_{k+1}\notin \mathcal G_{k}$ and thus are in $\mathcal G_{k+1}$.
To finish the argument, we need to find a  $\mathcal G_k\neq\emptyset$. This can be achieved by carefully analyzing the boundary points of $\mathcal N_j$. Beginning with an $\alpha_1$ we can see that it is either in $\mathcal G_1$ or $\Upsilon_1$. If it is in $\mathcal G_1$ we are set, so we assume that $|u(Z_k(\alpha_1),\alpha_1)|=\gamma_i$. Just below $\alpha_1$ solutions are in $\mathcal N_k$ for big $k$, and since $\bigcap_{n\in\mathbb{N}} \mathcal N_n = \emptyset$ we can find a $k_1$ and $\alpha^1_{k_1}$ that is an infimum of $(\alpha^1_{k_1},\alpha_1)\in \mathcal N_{k_1}$. This $\alpha^1_{k_1}$ must be in $\mathcal G_{k_1}$ or $\Upsilon_{k_1}$, if it is in $\Upsilon_{k_1}$ with  $|u(Z_k(\alpha^1_{k_1}),\alpha^1_{k_1})|=\gamma_j$, then $F(\gamma_j)<F(\alpha^1_{k_1})<F(\gamma_1)$, thus $j<i$. Repeating this argument through the decreasing maxima we can get an element in some $\mathcal G_{k_0}$.

\section{Multiplicity via abrupt magnitude changes of $f$}

In this section we want to study how abrupt changes in the size of $f$ generate multiplicity of solutions. We will approach this with functions that are defined by parts. This was addressed in \cite{cghh23} by the author in collaboration with  Carmen Cort\'azar and Marta Garc\'ia Huidobro, and we refer to this article for more details.

When $f_1$ satisfies some subcritical and superlinear conditions (for example $f_1=u^p-u$, $1<p< p^*$), solutions to the initial value problem are positive if the initial condition $u(0)$ is small, and change sign when $u(0)$ is big, with a ground state solution separating them at some initial value $u(0)=\alpha^*$.

\begin{figure}[h]
\begin{center}
\includegraphics[keepaspectratio, width=8cm]{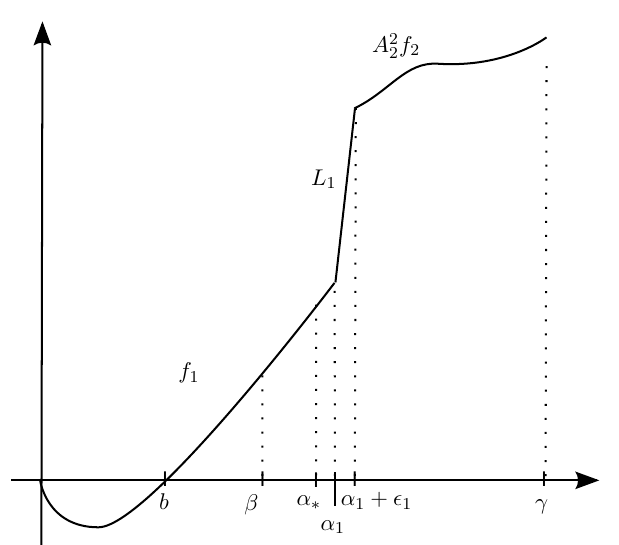}
 \end{center}
 \caption{The function $f$.}\label{Fig-fpp}
\end{figure}

  We consider the family of functions $f$ of the form
 \begin{eqnarray}\label{f}
f(s)=\begin{cases}
f_1(s) &  s\leq \alpha^*+\epsilon\\
L(s) & \alpha^*+\epsilon \leq s\leq \alpha^*+2\epsilon\\
A^2f_2(s) &  s\geq \alpha^*+2\epsilon
\end{cases}
\end{eqnarray}
where $f_2(s)$ is any positive continuous function and $L(s)$ is the line from $(\alpha^*+\epsilon, f_1(\alpha^*+\epsilon))$ to $(\alpha^*+2\epsilon, A^2f_2(\alpha^*+2\epsilon))$, and $A$, $\epsilon$ are suitable constants to be chosen. (See Figure \ref{Fig-fpp})

The problem with this function will have a ground state solution with initial condition $u(0)=\alpha^*$, and solutions that change sign when $\alpha^*<u(0)<\alpha^*+\epsilon$. In \cite{cghh23} we prove that, for $A$ big enough, there will be solutions with initial condition $u(0)>\alpha^*$ that will stay positive, an thus a second ground state solution.

We assume the following conditions on the nonlinearity $f$:
\begin{enumerate}
 \item[$(H_1)$]  $f_1\in C[0,\infty) \cap C^1(0,\infty)$, $f_1(0)=0$ and there exist $b\geq 0$ such that $f_1(s)>0$ for $s>b$,
	      $f_1(s)\leq 0$ for $s\in[0,b]$ and moreover
	      $f_1(s)<0$ on $(0,\epsilon)$ for some $\epsilon>0$;
also, by setting
$ F_1(s) = \int_0^s f_1(t) dt,
$
we assume that there exists a unique finite $\beta\geq b$ such that $F(\beta)=0$.
 \item[$(H_2)$] $(F_1/f_1)'(s) > (N-2)/(2N)$ for all $s>\beta$;
 \item[$(H_3)$] $\displaystyle f_1/(s-b)$ is increasing for all $s>b$.
\item[$(H_4)$] There is an initial condition $\alpha_*$ such that the problem
\begin{eqnarray}\label{pdef1}
\begin{gathered}
u''+\frac{N-1}{r}u'+f_1(u)=0,\quad r>0,\quad N> 2,\\
u(0)=\alpha_*,\quad u'(0)=0,
\end{gathered}
\end{eqnarray}
is a ground state solution.

\item[$(H_5)$] $f_2$ is a positive continuous function defined on $[\alpha_*,\gamma)$ for some $\alpha_*<\gamma\le\infty$.

\end{enumerate}

\medskip
Our main result in this section is the following.
\begin{theorem}\label{second solution}
Assume that $f_1$, $\alpha_*$ and $f_2$ satisfy the assumptions above.
Then, there exist positive constants $\bar\epsilon$ and  $\bar A$ such that for any $0<\epsilon_1<\bar\epsilon$ and $A_2>\bar A$, problem \eqref{pde} with $f$ given by \eqref{f} has at least two ground state solutions.
\end{theorem}

The second solution is found by proving that for an $\alpha_2>\alpha_*$ close to $\alpha_*$, we can choose small $\epsilon$ and large $A$ such that the solution with initial condition $u(0)=\alpha_2$ is positive for all $r$ (See Proposition \ref{EnP-f1} and what follows for more details). Thus there must be a ground state solution between $\alpha_1=\alpha_*+\epsilon \in \mathcal N_1$ and $\alpha_2\in \mathcal P_1$.

Moreover, if $f$ satisfies the subcritical type condition at infinity $(GST)$  seen in Section $2.1$, we can use the result by Gazzola, Serrin and Tang \cite{gst} to prove that for large initial conditions solutions must change sign. Therefore, there must be a third solution to \eqref{pde} above $\alpha_2$.

\begin{theorem}\label{third solution}
Assume that $f_1$, $\alpha_*$ and $f_2$ are as in Theorem \ref{second solution} with $\gamma=\infty$, and let {$\epsilon_1$ and $A_2$} be as in its conclusion. If $f$ satisfies $(GST)$, then problem \eqref{pde} has {at least three} ground state solutions.
\end{theorem}

Additionally, if after the second solution we make another jump in the value of $f$, this time making it small, we will have solutions with initial condition $u(0)>\alpha_2$ that will change sign, an thus a third ground state solution. We now consider $f$ given by
 \begin{eqnarray}\label{f_3}
f(s)=\begin{cases}
f_1(s) &  s\leq \alpha_1\\
L_1(s) & \alpha_1 \leq s\leq \alpha_1+\epsilon_1\\
A_2^2f_2(s) &   \alpha_1+\epsilon_1\leq  s\leq \alpha_2\\
L_2(s) & \alpha_2 \leq s\leq \alpha_2+\epsilon_2\\
A_3^2f_3(s) &  s\geq \alpha_2+\epsilon_2
\end{cases}
\end{eqnarray}
where  $L_1(s)$ is the line from $(\alpha_1, f_1(\alpha_1))$ to $(\alpha_1+\epsilon_1, A_2^2f_2(\alpha_1+\epsilon_1))$,  $L_2(s)$ is the line from $(\alpha_2, A_2^2f_2(\alpha_2))$ to $(\alpha_2+\epsilon_2, A_3^2f_3(\alpha_2+\epsilon_2))$ and  $\alpha_1,\  \epsilon_1$ and  $A_2$ are constants that satisfy Theorem \ref{second solution}.

 \begin{theorem}\label{solution f_3}
Under  assumptions $(H_1)$-$(H_5)$ ,
there exist positive  constants $\epsilon_1$, $\epsilon_2$ ,  $A_2$ and $A_3$ such that  problem \eqref{pde} with $f$ given by \eqref{f_3} has at least three ground state solutions.
\end{theorem}

 This follows easily from Lemma \ref{LemmaGST}, since when $f$ is small near the initial condition $\alpha_3$, the solution will have small slope and will cut $\alpha_2$ (and thus $\alpha_*$) with $r>C(\alpha_*)$. Then by Lemma \ref{LemmaGST} we have $\alpha_3\in \mathcal N_1$ and there must be a ground state solution between $\alpha_3$ and $\alpha_2\in \mathcal P_1$.

\begin{figure}[h]
\begin{center}
 \includegraphics[keepaspectratio, width=8cm]{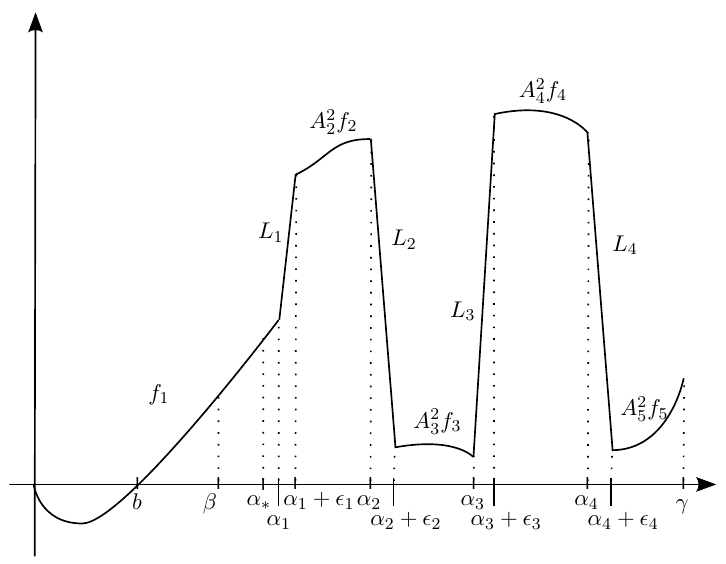}
 \end{center}
 \caption{The function $f$.}\label{Fig-fpp5}
 \end{figure}

 We can repeat this process making as many jumps as we want, obtaining as many ground state solutions as we want (see Figure \ref{Fig-fpp5}).

 \begin{theorem}\label{multiple}
Let   $f_i$, $i=1,\ldots,k$ satisfy assumptions $(H_1)$-$(H_5)$ For $i=2...k$, there exists constants  $\epsilon_i>0$, $ A_{i}>0$  and $\alpha_i$ with the condition $\alpha_*<\alpha_{i-1}+\epsilon_{i-1}<\alpha_{i}$  such that problem \eqref{pde} with
\begin{eqnarray*}
f=f_1\chi_{[0,\alpha_*+\epsilon_1]}+\sum_{i=2}^kL_{i-1}
\chi_{[\alpha_{i-1},\alpha_{i-1}+\epsilon_{i-1}]}+
\sum_{i=2}^{k-1}A_{i}^2f_{i}\chi_{[\alpha_{i-1}+\epsilon_{i-1},\alpha_{i}]}\\
+A_{k}^2f_{k}\chi_{[\alpha_{k-1}+\epsilon_{k-1},\gamma)}\end{eqnarray*}
has at least $k$ ground state solutions.
\end{theorem}

\begin{figure}[h]
\begin{center}
  \includegraphics[keepaspectratio,width=10cm]{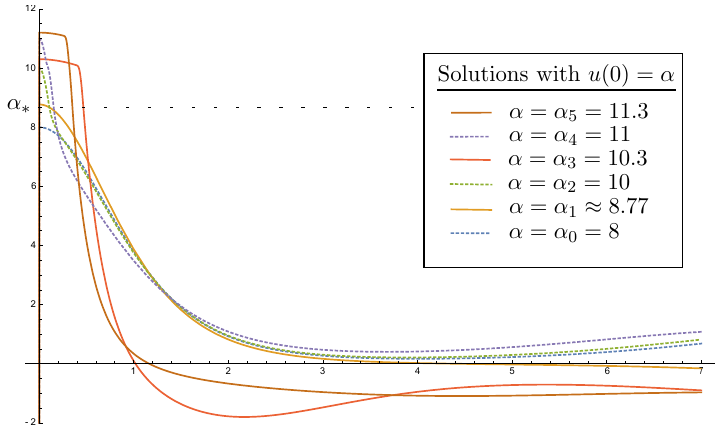}
  \end{center}
  \caption{Solutions of the example.}\label{EjA5}
\end{figure}

We can see the behavior described in these results on a numerical example.
\begin{example}
 In the special case when $N=4$, $f_1=u^2-u$, we estimate $\alpha_* \approx 8.672$ and choose $f_i=u^2$ for $i=2,\dots,5$, $\epsilon_i=\frac{1}{10}$, $A_2^2=A_4^2=10$, $A_3^2=A_5^2=\frac{1}{10}$ , (see  Figure \ref{EjA5}). We can see the different behavior of the solutions with $\alpha_i$ with even or odd $i$, showing the existence of ground state solutions between $\alpha_{i-1}$ and $\alpha_i$ for $i=1,2,3,4,5$.
\end{example}

The main ingredient of the proof of Theorem \ref{second solution} is the following proposition, that states that if a solution $u$ reaches the height $\alpha_*$ with small $r$ and controlled slope, it will be a positive solution with $u(0)\in \mathcal P_1$.

\begin{proposition}\label{EnP-f1}  Let $f_1$ satisfying the properties $(H_1)$-$(H_4)$,   and let $u$ be solution to
$$u''+\frac{N-1}{r}u'+f(u)=0$$
 that reaches the value $\alpha_*$ in $(H_4)$ at some $r_*\in(0,R)$.
Then given $\bar a,\ \bar b\in\mathbb R^+$, with  $ \bar b<(\alpha_*-b) (N-2) $, if $R>0$ is sufficiently small and $\bar a\le r_*|u'(r_*)|\le \bar b$, it holds that $u(r)>0$ for all  $r\ge r_*$ and $u(0)\in \mathcal P_1$.

\end{proposition}

\begin{figure}[h]
\begin{center}
 \includegraphics[keepaspectratio, width=8cm]{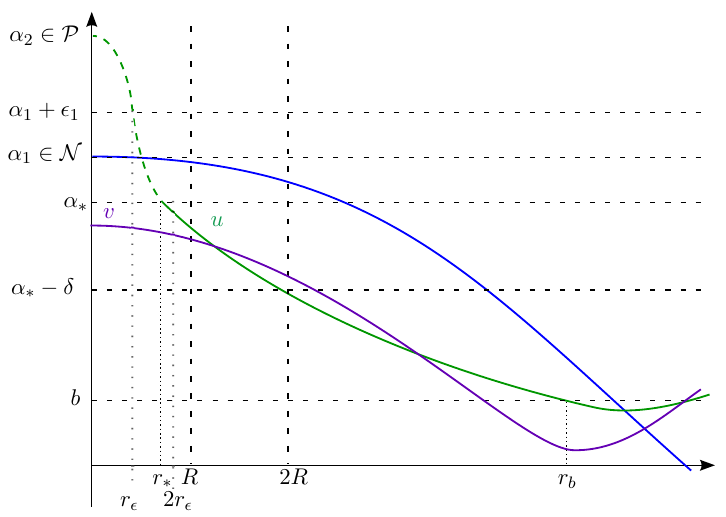}
 \end{center}
 \caption{ the functions $u$, $v$ and $w$.}\label{Fig-dem}
 \end{figure}

This is proved in several steps. We consider solutions $v$ of \eqref{ivp} with $v(0)$ close but below $\alpha_*$, thus $v(0)\in \mathcal P_1$ (see Figure \ref{Fig-dem}). It can be seen that these solution must intersect $u$,  and using the bounds on $r*$ and $r_*|u'(r_*)|$ as well as the super linear condition $(H3)$ to compare $u$ with two solutions $v$ and $w$ as above, we can find one such solution $v$ that intersects $u$ twice at values $u(r_I)>b$. Using the $P$ functional \eqref{PdeET} introduced by Erbe and Tang to compare $u$ and $v$ after the second intersection, and condition $(H2)$, we obtain that they can't intersect again before $Z_1$, and thus $u>v>0$ until it reaches a local minimum. Therefore $u\in \mathcal P_1$.

The proof of Theorem \ref{second solution} in \cite{cghh23} is concluded by observing that for any solution {$u$}, the function $ {v}(r)= {u}(Ar)$ satisfies:
$$v''(r)+\frac{(N-1)}{r}v'(r)= {A^2}\left(u''(Ar)+\frac{(N-1)}{Ar}u'(Ar)\right)= -{A^2}f(v(r))$$
 Moreover, $r_\epsilon^u= A r_\epsilon^v$ and $u'(r_\epsilon^u)= u'(Ar_\epsilon^v) =  v'(r_\epsilon^v)/A$ hence $r_\epsilon^u|u'(r_\epsilon^u)|$ is independent of $A$.

We begin with a solution $w$ of \eqref{ivp} for $f=f_2$, with $w(0)=\alpha_2$ close enough to $\alpha_*$ such that it has $ r_*|w'(r_*)|<(\alpha_*-b)(N-2) $. By choosing a big enough $A$, the solution with $f=A^2 f_2$ will have $r_*<R$ without changing $r_*|w'(r_*)|$, and thus the conditions for Proposition \ref{EnP-f1} are satisfied. Therefore, the solution of \eqref{ivp} for $f$ as in the Theorem with initial condition $u(0)=\alpha_2$ will be positive.

With this we conclude the proof of Theorem \ref{second solution}, and the review of recent results. It remains to study whether the magnitude changes generate $k$-bound state solutions, as well as how the place where the change occurs ($\alpha_*$) influence the result.




\end{document}